\documentclass{article}%
\usepackage{amsmath}%
\usepackage{amsfonts}%
\usepackage{amssymb}%
\usepackage{graphicx}
\usepackage[T1]{fontenc}
\usepackage[latin1]{inputenc}
\usepackage[french]{babel}

\newtheorem{theorem}{Théorème}

\newtheorem{definition}{Définition}

\newtheorem{lemma}{Lemme}

\newenvironment{proof}[1][Proof]{\textbf{#1.} }{\ \rule{0.5em}{0.5em}}

\begin{document}

\title{\'Eléments de distorsion du groupe des difféomorphismes isotopes à l'identité d'une variété compacte}
\author{E. Militon}
\date{\today}
\maketitle

\begin{abstract}
Dans cet article, on montre que, dans le groupe $\mathrm{Diff}_{0}^{\infty}(M)$ des difféomorphismes isotopes à l'identité d'une variété compacte $M$, tout élément récurrent est de distorsion. Pour ce faire, on généralise une méthode de démonstration utilisée par Avila pour le cas de $\mathrm{Diff}^{\infty}_{0}(\mathbb{S}^{1})$. La méthode nous permet de retrouver un résultat de Calegari et Freedman selon lequel tout homéomorphisme de la sphère isotope à l'identité est un élément de distorsion.
\end{abstract}

\section{\'Enoncé des résultats}

L'étude des éléments de distorsion des groupes de difféomorphismes ou d'homéomorphismes non-conservatifs d'une variété trouve son origine dans une question de Franks et Handel (voir $\cite{FH}$ et $\cite{F}$) : une rotation du cercle est-elle distordue dans le groupe des homéomorphismes du cercle ou des difféomorphismes du cercle ? Qu'en est-il si l'on remplace le cercle par la sphère de dimension 2 ?

La réponse à ces questions est fournie par Calegari et Freedman dans $\cite{CF}$ : une rotation du cercle est distordue dans le groupe des difféomorphismes de classe $C^{1}$ du cercle. Cependant, les auteurs précisent ne pas savoir s'il en est de même en régularité $C^{\infty}$. Dans le même article, ils prouvent qu'une rotation de la sphère $S^{2}$ est distordue dans le groupe des difféomorphismes de classe $C^{\infty}$ de la sphère. Enfin, en ce qui concerne la régularité $C^{0}$, Calegari et Freedman ont prouvé un résultat très général : en toute dimension $N$, un homéomorphisme $h$ de la sphère $S^{N}$ est distordu.

Dans un article ultérieur, Avila a montré que, dans le groupe des difféomorphismes du cercle de classe $C^{\infty}$, tout élément récurrent est distordu. Nous généralisons dans cet article le résultat d'Avila à toute variété. 

\section{Résultats}

Avant toute chose, commençons par introduire des définitions et des notations qui nous seront utiles par la suite.

\begin{definition}
Soit $G$ un groupe. Pour une partie finie $S$ de $G$, si un élément $g$ de $G$ est dans le groupe engendré par $S$, on note :
$$l_{S}(g)=inf \left\{n \in \mathbb{N}, \ \exists (\epsilon_{i})_{1 \leq i \leq n} \in \left\{\pm 1 \right\}^{n}, \exists (s_{i})_{1 \leq i \leq n} \in S^{n}, \ g = \prod_{i=1}^{n}s_{i}^{\epsilon_{i}} \right\}.$$
Un élément $g$ de $G$ est dit distordu (ou de distorsion) s'il existe une partie finie $S$ de $G$ telle que $g$ appartient au groupe engendré par $S$ et :
$$\lim_{n \rightarrow +\infty} \frac{l_{S}(g^{n})}{n}=0.$$
\end{definition}

On remarque que, comme la suite $(l_{S}(g^{n}))_{n \in \mathbb{N}}$ est sous-additive, il suffit de montrer que :
$$ \liminf_{n \rightarrow + \infty} \frac{l_{S}(g^{n})}{n}=0$$
pour obtenir que l'élément $g$ est distordu.

Etant donnée une variété différentiable $M$, on note :
\begin{itemize}
\item $\mathrm{Homeo}_{0}(M)$ l'ensemble des homéomorphismes à support compact dans $M$ isotopes à l'identité par une isotopie à support compact;
\item pour un élément $r$ de $\mathbb{N} \cup \left\{ \infty \right\}$, $\mathrm{Diff}_{0}^{r}(M)$ l'ensemble des difféomorphismes de $M$ de classe $C^{r}$ isotopes à l'identité par une isotopie à support compact (en particulier, $\mathrm{Homeo}_{0}(M)=\mathrm{Diff}_{0}^{0}(M)$).
\end{itemize}

Le support d'un homéomorphisme $f$ de $M$ est ici défini par :
$$supp(f)=\overline{ \left\{ x \in M, \ f(x) \neq x \right\} }.$$
\\
Si $f$ et $g$ sont deux éléments de $\mathrm{Diff}_{0}^{r}(M)$, on note $d_{r}$ une distance qui définit la topologie de $\mathrm{Diff}_{0}^{r}(M)$.

\begin{definition}
On dit qu'un difféomorphisme $f$ de $\mathrm{Diff}_{0}^{\infty}(M)$ est récurrent si et seulement si :
$$ \liminf_{n \rightarrow + \infty} d_{\infty}(f^{n}, Id_{M})=0.$$
\end{definition}

L'objet de la présente note est de démontrer les théorèmes suivants. Ce premier théorème généralise le résultat d'Avila pour les difféomorphismes du cercle (voir $\cite{Av}$).

\begin{theorem}
Si $M$ est une variété compacte, tout élément récurrent de \\ $\mathrm{Diff}_{0}^{\infty}(M)$ est distordu. 
\end{theorem}

Comme dans $\cite{Av}$, la méthode employée permet de démontrer que toute suite d'éléments récurrents $(f_{n})_{n \in \mathbb{N}}$ de $\mathrm{Diff}_{0}^{\infty}(M)$ est simultanément distordue, au sens où l'on peut trouver un ensemble fini $S$ tel que :
$$ \forall n \in \mathbb{N}, \ \liminf_{ p \rightarrow + \infty} \frac{l_{S}(f_{n}^{p})}{p}=0.$$
La méthode employée par $\cite{Av}$ permet de donner une nouvelle démonstration du résultat de Calegari et Freedman (voir $\cite{CF}$).

\begin{theorem}
(Calegari-Freedman) Tout élément de $\mathrm{Homeo}_{0}(\mathbb{S}^{n})$ est distordu.
\end{theorem}

Là encore, on pourrait montrer que toute suite d'éléments de $\mathrm{Homeo}_{0}(\mathbb{S}^{n})$ est simultanément distordue.

\section{Démonstration des théorèmes}

La démonstration des théorèmes 1 et 2 repose sur une généralisation des résultats de $\cite{Av}$ à toute variété. La démarche et les notations sont similaires à celles présentées dans $\cite{Av}$ mais la généralisation de la méthode aux dimensions supérieures nécessite de manière cruciale des résultats de perfection locale qui proviennent de $\cite{HT}$. On ne démontrera le théorème 1 que dans le cas d'une variété $M$ de dimension supérieure ou égale à 2. Le cas de la dimension 1 est traité dans $\cite{Av}$.\\
Les théorèmes vont découler des lemmes suivants.

\begin{lemma}
Soit $M$ une variété compacte de dimension supérieure ou égale à $2$. \\
Il existe des suites $(\epsilon_{n})_{n \in \mathbb{N}}$ et $(k_{n})_{n \in \mathbb{N}}$ de réels strictement positifs telles que, toute suite $(h_{n})_{n \in \mathbb{N}}$ d'éléments de $\mathrm{Diff}_{0}^{\infty}(M)$ telle que : 
$$\forall n \in \mathbb{N}, \ d_{\infty}(h_{n},Id_{M})<\epsilon_{n},$$
vérifie la propriété suivante : il existe un ensemble fini $S \subset \mathrm{Diff}_{0}^{\infty}(M)$ tel que, pour tout entier $n$,
\begin{itemize}
\item le difféomorphisme $h_{n}$ appartient au sous-groupe engendré par $S$.
\item on a l'inégalité : $l_{S}(h_{n}) \leq k_{n}$.
\end{itemize}
\end{lemma}

\begin{lemma}
Soit $N$ un entier naturel.\\
Il existe une suite $(k_{n})_{n \in \mathbb{N}}$ de réels positifs telle que, pour toute suite $(h_{n})_{n \in \mathbb{N}}$ d'éléments de $\mathrm{Homeo}_{0}(\mathbb{S}^{N})$, il existe un ensemble fini $S \subset \mathrm{Homeo}_{0}(\mathbb{S}^{N})$ tel que, pour tout entier n :
\begin{itemize}
\item l'homéomorphisme $h_{n}$ appartient au groupe engendré par $S$.
\item on a l'inégalité : $l_{S}(h_{n}) \leq k_{n}$.
\end{itemize}
\end{lemma}

Admettons pour l'instant ces lemmes, qui seront démontrés dans la section suivante, et démontrons les théorèmes.

\noindent \begin{proof}[Démonstration des théorèmes] Soit $f$ un élément récurrent de $\mathrm{Diff}_{0}^{\infty}(M)$. Considérons une application strictement croissante $p: \mathbb{N} \rightarrow \mathbb{N}$ telle que :
$$ \left\{
\begin{array}{l}
\lim_{n \rightarrow +\infty}\frac{k_{n}}{p(n)}=0 \\
d_{\infty}(f^{p(n)},Id_{M}) \leq \epsilon_{n}
\end{array}
\right. ,
$$
où les suites $(\epsilon_{n})_{n \in \mathbb{N}}$ et $(k_{n})_{n \in \mathbb{N}}$ sont données par le lemme 1. Ce même lemme appliqué à la suite $(f^{p(n)})_{n \in \mathbb{N}}$ nous donne l'existence d'un ensemble fini $S$ qui montre que $f$ est distordu. Le théorème 2 se montre de la même manière, en utilisant le lemme 2.
\end{proof}

\section{Démonstration des lemmes 1 et 2}

Pour mener à bien la démonstration de ces lemmes, nous aurons besoin des lemmes suivants, qui seront démontrés dans la section suivante et portent sur les difféomorphismes de $\mathbb{R}^{n}$. Les lemmes 3 et 4 sont des analogues des lemmes 1 et 2 dans le cas des commutateurs de difféomorphismes de $\mathbb{R}^{N}$. Notons $B(0,2)$ la boule ouverte de centre $0$ et de rayon $2$. Si $f$ et $g$ sont des difféomorphismes de $\mathbb{R}^{n}$, on note $[f,g]$ le difféomorphisme $fgf^{-1}g^{-1}$.

\begin{lemma}
Soient $N$ un entier naturel et $r$ un élément de $\mathbb{N} \cup \left\{\infty \right\}$.\\
Il existe des suites $(\epsilon'_{n})_{n \in \mathbb{N}}$ et $(k'_{n})_{n \in \mathbb{N}}$ de réels strictement positifs vérifiant la propriété suivante. Soient $(f_{n})_{n \in \mathbb{N}}$ et $(g_{n})_{n \in \mathbb{N}}$ des suites d'éléments de $\mathrm{Diff}_{0}^{r}(\mathbb{R}^{N})$ à support inclus dans la boule $B(0,2)$ telles que :
$$\forall n \in \mathbb{N}, \ d_{r}([f_{n},g_{n}],Id_{\mathbb{R}^{N}})<\epsilon'_{n}.$$
Alors il existe un ensemble fini $S'$ inclus dans $\mathrm{Diff}_{0}^{r}(\mathbb{R}^{N})$ tel que pour tout entier $n$, le difféomorphisme $[f_{n},g_{n}]$ appartient au groupe engendré par $S'$ et $l_{S}([f_{n},g_{n}]) \leq k'_{n}$.
\end{lemma}

Dans le cas particulier où $r=0$, on a un lemme un peu plus fort.

\begin{lemma}
Il existe une suite $(k'_{n})_{n \in \mathbb{N}}$ de réels strictement positifs telle que, pour toute suite $(\tilde{h}_{n})_{n \in \mathbb{N}}$ d'éléments de $\mathrm{Homeo}_{0}(\mathbb{R}^{N})$ à support dans $B(0,2)$, il existe un ensemble fini $S'$ inclus dans $\mathrm{Homeo}_{0}(\mathbb{R}^{N})$ tel que, pour tout entier $n$, l'homéomorphisme $\tilde{h}_{n}$ appartient au groupe engendré par $S'$ et $l_{S}(\tilde{h}_{n}) \leq k'_{n}$.
\end{lemma}

\bigskip

\noindent \begin{proof}[Démonstration du lemme 1] On note $N$ la dimension de $M$ ($N \geq 2$). On considère un recouvrement ouvert $(U_{i})_{0 \leq i \leq p}$ de $M$ constitué d'ouverts difféomorphes à $\mathbb{R}^{N}$ dont l'adhérence est incluse dans un ouvert de carte de $M$. Pour tout entier $i$ entre $0$ et $p$, on choisit $\varphi_{i}$ une carte de $M$ définie sur un voisinage de l'adhérence de $U_{i}$ qui vérifie : 
$$\varphi_{i}(U_{i}) \subset B(0,2).$$
Notons $\psi$ une bijection de $\mathbb{N} \times \left\{1,...,p\right\}\times \left\{1,\ldots, 4N\right\}$ sur $\mathbb{N}$.

On va maintenant construire la suite $(\epsilon_{n})_{n \in \mathbb{N}}$ recherchée à l'aide de la suite $(\epsilon'_{l})_{l \in \mathbb{N}}$ donnée par le lemme 3 dans le cas $r=+\infty$.\\
D'après le lemme de perfection lisse (voir appendice), pour tout entier naturel $n$, on peut choisir $\epsilon_{n}$ suffisamment petit de sorte que, si un difféomorphisme $h$ de $\mathrm{Diff}_{0}^{\infty}(M)$ vérifie $d(h,Id_{M}) < \epsilon_{n}$, alors il existe deux familles $(f_{k,l})_{0 \leq k \leq p,1 \leq l \leq 3N}$ et $(g_{k,l})_{0 \leq k \leq p,1 \leq l \leq 3N}$ d'éléments de $\mathrm{Diff}_{0}^{\infty}(M)$, où $g_{k,l}$ et $f_{k,l}$ sont supportés dans l'ouvert $U_{k}$ telles que :
$$ \left\{
\begin{array}{l}
h=\prod_{k=0}^{p}\prod_{l=1}^{3N}[f_{k,l},g_{k,l}] \\
\forall k,l \in [0,p]\times[1,3N], \ d_{\infty}(\varphi_{k} \circ f_{k,l} \circ \varphi_{k}^{-1},Id_{\mathbb{R}^{N}})< \epsilon'_{\psi(n,k,l)} \\
\forall k,l \in [0,p]\times[1,3N], \ d_{\infty}(\varphi_{k} \circ g_{k,l} \circ \varphi_{k}^{-1},Id_{\mathbb{R}^{N}})< \epsilon'_{\psi(n,k,l)}
\end{array}
\right.
.
$$
Donnons-nous une suite $(h_{n})_{n \in \mathbb{N}}$ d'éléments de $\mathrm{Diff}_{0}^{\infty}(M)$ qui vérifie :
$$\forall n \in \mathbb{N}, \ d_{\infty}(h_{n}, Id_{M}) < \epsilon_{n}.$$
La définition des $\epsilon_{n}$ nous donne deux familles $(\tilde{f}_{n,k,l})$ et $(\tilde{g}_{n,k,l})$ associées à $h_{n}$. On pose alors :
$$f_{\psi(n,k,l)}=\varphi_{k} \circ \tilde{f}_{n,k,l} \circ \varphi_{k}^{-1}$$
$$g_{\psi(n,k,l)}=\varphi_{k} \circ \tilde{g}_{n,k,l} \circ \varphi_{k}^{-1}$$
et on applique le lemme 3 aux suites $(f_{n})_{n \in \mathbb{N}}$ et $(g_{n})_{n \in \mathbb{N}}$ obtenues, ce qui achève la démonstration du lemme 1 en posant :
$$k_{n}= \sum_{k,l} k'_{\psi(n,k,l)}.$$
\end{proof}

\bigskip

\noindent \begin{proof}[Démonstration du lemme 2] On utilise le lemme suivant qui découle d'un résultat très profond dû à Kirby, Siebenmann et Quinn pour le cas de la dimension supérieure ou égale à $4$ (voir $\cite{CF}$ lemme 6.10, $\cite{Kir}$ et $\cite{Q}$) :

\begin{lemma}
On considère deux disques fermés inclus dans $S^{n}$ dont les intérieurs recouvrent $S^{n}$. Alors tout homéomorphisme $h$ isotope à l'identité s'écrit comme produit de six éléments de $\mathrm{Homeo}_{0}(S^{n})$ qui sont chacun à support inclus dans l'un de ces disques.    
\end{lemma}

En utilisant ce résultat, la même méthode que précédemment, conjuguée au lemme 4, donne le lemme 2.
\end{proof}

\section{Démonstration des lemmes 3 et 4}
Là encore, elle suit la même méthode que celle présentée dans $\cite{Av}$.

\bigskip

\noindent \begin{proof}[Démonstration du lemme 3.] On note $F_{1}$ un élément de $\mathrm{Diff}_{0}^{\infty}(\mathbb{R}^{N})$ qui vérifie :
$$ \forall x \in B(0,2), \ F_{1}(x)= \lambda x ,$$
avec $0<\lambda<1$.\\
Soit $F_{2}$ un élément de $\mathrm{Diff}_{0}^{\infty}(\mathbb{R}^{N})$ à support inclus dans $B(0,2)$ qui vérifie :
$$ \forall x \in B(0,1), F_{2}(x)=x+a,$$
où $a$ appartient à $\mathbb{R}^{N}-\left\{0\right\}$ et est de norme strictement inférieure à $1$.\\
Soit $F_{3}$ un élément de $\mathrm{Diff}_{0}^{\infty}(\mathbb{R}^{N})$ à support inclus dans $B(0,2)$ qui vérifie :
\begin{itemize}
\item la suite $(F_{3}^{n}(0))_{n \in \mathbb{N}}$ est constituée d'éléments deux à deux distincts.
\item la suite $(F_{3}^{n}(0))_{n \in \mathbb{N}}$ converge vers $x_{0}=(1,0, \ldots, 0)$.
\end{itemize}
Considérons une suite d'entiers $(l_{n})_{n \in \mathbb{N}}$ croissant suffisamment vite pour que :
\begin{itemize}
\item $n \neq n' \Rightarrow F_{3}^{n}F_{1}^{l_{n}}(B(0,2)) \cap F_{3}^{n'}F_{1}^{l_{n'}}(B(0,2))= \emptyset$.
\item le diamètre de $F_{3}^{n}F_{1}^{l_{n}}(B(0,2))$ converge vers 0.
\end{itemize}
On note $F_{n}=F_{3}^{n}F_{1}^{l_{n}}$, $U_{n}=F_{3}^{n}F_{1}^{l_{n}}(B(0,2))$ et $\hat{F_{n}}=F_{n}F_{2}F_{n}^{-1}$.\\
On considère pour chaque entier $n$ un ouvert $V_{n}$ qui vérifie :
\begin{itemize}
\item $F_{3}^{n}(0) \in V_{n} \subset U_{n}$.
\item $\hat{F_{n}}(V_{n}) \cap V_{n}= \emptyset$.
\end{itemize}
On considère aussi une suite d'entiers $(\tilde{l}_{n})_{n \in \mathbb{N}}$ de sorte que :
$$F_{3}^{n}F_{1}^{\tilde{l}_{n}}(B(0,2)) \subset V_{n}$$
et on note $\tilde{F}_{n}= F_{3}^{n}F_{1}^{\tilde{l}_{n}}$.

Donnons-nous des suites $(f_{n})_{n \in \mathbb{N}}$ et $(g_{n})_{n \in \mathbb{N}}$ d'éléments de $\mathrm{Diff}_{0}^{r}(\mathbb{R}^{N})$ à support dans $B(0,2)$ et choisissons la suite $(\epsilon'_{n})_{n \in \mathbb{N}}$ tendant suffisamment vite vers $0$ de sorte que, si $d_{r}(f_{n}, Id) < \epsilon'_{n}$ et $d_{r}(g_{n},Id)< \epsilon'_{n}$ pour tout entier $n$, alors les applications $F_{4}, F_{5} : \mathbb{R}^{N} \rightarrow \mathbb{R}^{N}$ définies par :
$$\begin{array}{c}
\forall n \in \mathbb{N}, \ F_{4 | V_{n}}=\tilde{F}_{n}f_{n}\tilde{F}_{n}^{-1} \\
\forall n \in \mathbb{N}, \ F_{5 | V_{n}}=\tilde{F}_{n}g_{n}\tilde{F}_{n}^{-1} \\
F_{4 | \mathbb{R}^{N}- \bigcup_{n \in \mathbb{N}}V_{n}}=F_{5 | \mathbb{R}^{N}- \bigcup_{n \in \mathbb{N}}V_{n}}=Id_{\mathbb{R}^{N}- \bigcup_{n \in \mathbb{N}}V_{n}}
\end{array}
$$
sont de classe $C^{r}$ en $x_{0}$ (et ce sont alors des $C^{r}$-difféomorphismes de $M$ à support inclus dans $B(0,2)$ isotopes à l'identité).

Remarquons que
$$A_{n}=F_{4}\hat{F}_{n}F_{4}^{-1}\hat{F}_{n}^{-1}$$
est à support dans $V_{n} \cup \hat{F}_{n}(V_{n})$, vaut $\tilde{F}_{n}f_{n}\tilde{F}_{n}^{-1}$ sur $V_{n}$ et $\hat{F}_{n}\tilde{F}_{n}f_{n}^{-1}\tilde{F}_{n}^{-1}\hat{F}_{n}^{-1}$ sur $\hat{F}_{n}(V_{n})$. On peut faire une remarque analogue pour les applications $B_{n}=F_{5}\hat{F}_{n}F_{5}^{-1}\hat{F}_{n}^{-1}$ et $C_{n}=F_{4}^{-1}F_{5}^{-1}\hat{F}_{n}F_{5}F_{4}\hat{F}_{n}^{-1}$. On a alors :
$$ A_{n}B_{n}C_{n}=\tilde{F}_{n}[f_{n},g_{n}] \tilde{F}_{n}^{-1}.$$
Ainsi :
$$[f_{n},g_{n}]=\tilde{F}_{n}^{-1}A_{n}B_{n}C_{n}\tilde{F}_{n}.$$
En prenant $S= \left\{F_{1}, F_{2}, F_{3}, F_{4}, F_{5} \right\}$, on obtient le résultat escompté.
\end{proof}

\bigskip 

\noindent \begin{proof}[Démonstration du lemme 4] Remarquons que, dans la démonstration précédente, l'apparition de la suite $(\epsilon_{n})_{n \in \mathbb{N}}$ est liée à un défaut de régularité de $F_{4}$ et de $F_{5}$. En régularité $C^{0}$, ce problème n'apparait pas, ce qui démontre le lemme 4 dans le cas où la suite $(\tilde{h}_{n})_{n \in \mathbb{N}}$ est constituée de commutateurs. Il suffit ensuite de démontrer le lemme suivant.
\end{proof}

\bigskip

\begin{lemma}
Tout homéomorphisme de $\mathbb{R}^{N}$ à support compact isotope à l'identité s'écrit comme un commutateur dans $\mathrm{Homeo}_{0}(\mathbb{R}^{N})$.
\end{lemma}

\noindent \begin{proof}[Démonstration (tirée de $\cite{Bo}$, démonstration du théorème 1.1.3)]  On note $\varphi$ un homéomorphisme de $\mathrm{Homeo}_{0}(\mathbb{R}^{n})$ de restriction à $B(0,2)$ définie par :
$$\begin{array}{rcl}
B(0,2) & \rightarrow & \mathbb{R}^{N} \\
x & \mapsto & \frac{x}{2}
\end{array}$$
Pour tout entier naturel $n$, on note :
$$ A_{n}=\left\{x \in \mathbb{R}^{N}, \frac{1}{2^{n+1}}\leq \left\|x\right\| \leq \frac{1}{2^{n}} \right\}.$$
Soit $h$, un élément de $\mathrm{Homeo}_{0}(\mathbb{R}^{N})$. Comme tout élément de $\mathrm{Homeo}_{0}(\mathbb{R}^{N})$ est conjugué à un élément à support inclus dans l'intérieur de $A_{0}$, on peut supposer $h$ à support inclus dans l'intérieur de $A_{0}$. On définit alors $g \in \mathrm{Homeo}_{0}(\mathbb{R}^{N})$ par :
\begin{itemize}
\item $g=Id$ en dehors de $B(0,1)$.
\item pour tout entier naturel $i$, $g_{|A_{i}}=\varphi^{i}h_{i} \varphi^{-i}$.
\item $g(0)=0$.
\end{itemize}
Alors :
$$ h= [g,\varphi].$$ 
\end{proof}

\appendix

\section{Appendice : perfection locale de $\mathrm{Diff}_{0}^{\infty}(M)$}

L'objet de cet appendice est de démontrer le résultat suivant qui a été utilisé au cours de la démonstration du lemme 1 :

\begin{theorem}
Soit $M$, une variété compacte connexe de dimension $n$. Fixons un recouvrement ouvert $(U_{k})_{0 \leq k \leq p}$ de $M$ constitué d'ouverts d'adhérences difféomorphes à la boule unité de $\mathbb{R}^{n}$. Alors, pour tout voisinage $\Omega$ de l'identité dans $\mathrm{Diff}_{0}^{\infty}(M)$, il existe un voisinage $\Omega'$ de l'identité dans $\mathrm{Diff}_{0}^{\infty}(M)$ tel que, pour tout difféomorphisme $f$ de $\Omega'$, il existe des familles de difféomorphismes $(f_{k,l})_{0 \leq k \leq p, 1 \leq l \leq 3n}$ et $(g_{k,l})_{0 \leq k \leq p, 1 \leq l \leq 3n}$ dans $\Omega$ tels que :
$$ f = \prod_{k=0}^{p} \prod_{l=1}^{3n} [f_{k,l},g_{k,l}]$$
et les difféomorphismes $f_{k,l}$ et $g_{k,l}$ sont à support dans $U_{k}$.
\end{theorem}

Cette démonstration est une réalisation élémentaire de l'idée de Stefan Haller et Josef Teichmann de décomposer un difféomorphisme en produit de difféomorphismes préservant certains feuilletages (voir $\cite{HT}$). Elle repose de manière essentielle sur le théorème KAM d'Herman sur les difféomorphismes du cercle. On remarque que cette propriété démontre la perfection de $\mathrm{Diff}_{0}^{\infty}(M)$ et donc la simplicité de ce groupe. C'est donc une alternative à la preuve de Thurston et Mather (voir $\cite{Bo}$ ou $\cite{Ban}$). La démonstration donne aussi la perfection locale (et donc la perfection) de $\mathrm{Diff}_{0}^{\infty}(\mathbb{R}^{n})$ pour $n$ supérieur ou égal à 2 mais ne permet pas de conclure dans le cas $n=1$. 

D'après le lemme de fragmentation (voir $\cite{HT}$ proposition 1 ou $\cite{Bo}$ théorème 2.2.1 pour une démonstration), pour tout réel $\eta>0$, il existe $\alpha >0$ tel que, pour tout élément $h$ de $\mathrm{Diff}_{0}^{\infty}(M)$, si $d(h,Id_{M})< \alpha$, alors il existe une famille $(h_{i})_{0 \leq i \leq p}$ d'éléments de $\mathrm{Diff}_{0}^{\infty}(M)$ telle que :
$$\left\{
\begin{array}{l}
supp(h_{i}) \subset U_{i} \\
h = \prod_{i=0}^{p} h_{i} \\
\forall i \in [0,p]\cap \mathbb{N}, \ d(h_{i}, Id_{M}) < \eta
\end{array}
\right.
.$$

Il reste  à effectuer la construction suivante. On se donne $U$, $V$, deux ouverts de $\mathbb{R}^{n}$, où $U$ est un cube d'adhérence incluse dans $V$, et un voisinage $\Omega$ de l'identité dans $\mathrm{Diff}_{0}^{\infty}(V)$.
Nous allons montrer l'existence d'un voisinage $\Omega'$ de l'identité dans $\mathrm{Diff}_{0}^{\infty}(U)$ tel que, pour tout difféomorphisme $f$ de $\Omega'$, il existe des difféomorphismes $f_{1}, f_{2}, \ldots, f_{3n}$ dans $\Omega$ tels que :
$$ f = [f_{1},f_{2}] \circ [f_{3},f_{4}] \circ \ldots \circ [f_{3n-1}, f_{3n}].$$

Pour montrer cette dernière propriété, la stratégie sera la suivante : on va commencer par décomposer un difféomorphisme $f$ proche de l'identité en tant que produit de $n$ difféomorphismes qui préservent chacun les feuilles d'un feuilletage en droites. Chacun de ces feuilletages en droites de $U$ va être considéré comme une partie d'un feuilletage en cercles d'un anneau inclus dans $V$. Il suffira ensuite d'appliquer (soigneusement) le théorème d'Herman sur les difféomorphismes du cercle pour conclure que chacun des difféomorphismes apparaissant dans la décomposition de $f$ s'écrit comme produit de deux commutateurs constitués d'éléments qui peuvent être choisis aussi proche que l'on veut de l'identité tant que $f$ est suffisamment proche de l'identité.

Détaillons maintenant les arguments ci-dessus. Pour un entier $k$ entre $1$ et $n$, on note $F_{k}$ le feuilletage constitué de l'ensemble des droites parallèles au k-ième axe de coordonnées. On note $\mathrm{Diff}_{0}^{\infty}(U, F_{k})$ l'ensemble des difféomorphismes de $U$ à support compact et compactement isotopes à l'identité qui préservent les feuilles du feuilletage $F_{k}$. Construisons par récurrence sur $k$ une application définie et continue sur un voisinage de l'identité $\phi_{k} : \mathrm{Diff}_{0}^{\infty}(U) \rightarrow \mathrm{Diff}_{0}^{\infty}(U,F_{k})$ telle que $\phi_{k}(Id_{U})=Id_{U}$ et telle que, pour un difféomorphisme $f$ suffisamment proche de l'identité :  
$$p_{k} \circ f \circ \phi_{1}(f)^{-1} \circ \phi_{2}(f)^{-1} \circ \ldots \circ \phi_{k}(f)^{-1}=p_{k},$$
où $p_{k}$ désigne la projection sur les $k$ premières coordonnées. Supposons $\phi_{1},\phi_{2}, \ldots, \phi_{k}$ construites. Posons, pour un difféomorphisme $f$ de $\mathrm{Diff}_{0}^{\infty}(U)$ proche de l'identité et pour $x$ dans $u$:
$$f \circ \phi_{1}(f)^{-1} \circ \phi_{2}(f)^{-1} \circ \ldots \circ \phi_{k}(f)^{-1}(x)=(x_{1}, \ldots, x_{k}, f_{k+1}(x), f_{k+2}(x), \ldots, f_{n}(x)).$$
On définit alors $\phi_{k+1}$ par :
$$ \forall x \in U, \phi_{k+1}(f)(x)=(x_{1}, \ldots, x_{k}, f_{k+1}(x), x_{k+2}, \ldots, x_{n})$$
pour $f$ suffisamment proche de l'identité, ce qui conclut la récurrence. Les applications $\phi_{k}$ ainsi construites sont $C^{\infty}$, valent l'identité en l'identité et vérifient :
$$f=\phi_{n}(f) \circ \phi_{n-1}(f) \circ \ldots \circ \phi_{1}(f).$$

Fixons un entier $k$. On considère comme prévu un plongement :
$$ \psi_{k} : \mathbb{R}^{k-1} \times \mathbb{R}/\mathbb{Z} \times \mathbb{R}^{n-k} \rightarrow V$$
qui envoie $(-1,1)^{k-1} \times (-\frac{1}{4},\frac{1}{4}) \times (-1,1)^{n-k}$ sur $U$ et le feuilletage constitué des droites parallèles au k-ième axe de coordonnées de $(-1,1)^{k-1} \times (-\frac{1}{4},\frac{1}{4}) \times (-1,1)^{n-k}$ sur le feuilletage constitué des droites parallèles au k-ième axe de coordonnées de $U$. 
Le difféomorphisme $\psi_{k}^{-1} \circ \phi_{k}(f) \circ \psi_{k}$ est alors identifiable à une famille (à $n-1$ paramètres) de difféomorphismes du cercle. On écrit pour des éléments $(x,t,y)$ de $\mathbb{R}^{k-1} \times \mathbb{R} /\mathbb{Z} \times \mathbb{R}^{n-k}$ :
$$\psi_{k}^{-1} \circ \phi_{k}(f) \circ \psi_{k}(x,t,y)=(x,g(x,y)(t),y),$$
où $g:\mathbb{R}^{n-1} \rightarrow \mathrm{Diff}_{0}^{\infty}(\mathbb{R}/\mathbb{Z})$ est une application continue qui vérifie :
$$g(x,y)=Id$$
lorsque $(x,y) \notin [-1,1]^{n-1}$.

Il reste à écrire $g(x,y)$ en tant que produit de commutateurs d'éléments de $\mathrm{Diff}_{0}^{\infty}(\mathbb{R}/\mathbb{Z})$ proches de l'identité qui dépendent de manière $C^{\infty}$ de $x$ et de $y$, dont les dérivées par rapport à $x$ et à $y$ sont petites et qui valent l'identité lorsque $(x,y) \notin (-2,2)^{n-1}$.

Le lemme suivant est dû à Haller et Teichmann (voir $\cite{HT2}$, exemple 3) :

\begin{lemma} 
Pour tout voisinage W' l'identité dans $\mathrm{Diff}_{0}^{\infty}(\mathbb{R}/\mathbb{Z})$, il existe un voisinage $W$ de l'identité dans $\mathrm{Diff}_{0}^{\infty}(\mathbb{R}/\mathbb{Z})$, des difféomorphismes $A$, $B$ et $C$ de $W'$ et des applications 
$$a,b,c : W \rightarrow \mathrm{Diff}_{0}^{\infty}(\mathbb{R}/\mathbb{Z})$$
de classe $C^{\infty}$, qui valent l'identité en l'identité, et vérifient, pour tout $h$ de $W$ :
$$h=[a(h),A][b(h),B][c(h),C].$$ 
\end{lemma}

\paragraph{Remarque.} La démonstration de ce résultat repose sur le théorème d'Hermann pour les difféomorphismes du cercle. De plus, on peut supposer que $A$, $B$, $C$ sont des éléments de $PSL_{2}(\mathbb{R}) \subset \mathrm{Diff}_{0}^{\infty}(\mathbb{R}/\mathbb{Z})$. 

\paragraph{Fin de la démonstration du théorème.} Notons $\lambda : \mathbb{R}^{n-1} \rightarrow \mathbb{R}$ une application $C^{\infty}$ à support dans $(-2,2)^{n-1}$ qui vaut $1$ sur un voisinage de $[-1,1]^{n-1}$. On choisit alors un voisinage $W'$ de l'identité dans $\mathrm{Diff}_{0}^{\infty}(\mathbb{R}/\mathbb{Z})$ de sorte que, pour tout difféomorphisme $D$ de $W'$, si l'on note $\tilde{D}$ l'application définie par :
$$\forall (x,t,y) \in \mathbb{R}^{k-1} \times \mathbb{R}/\mathbb{Z} \times \mathbb{R}^{n-k}, \ \tilde{D}(x,t,y)=(x,\lambda(x,y)D(t)+(1-\lambda(x,y))t,y),$$
alors $\tilde{D}$ est un difféomorphisme de $\mathbb{R}^{k-1} \times \mathbb{R}/\mathbb{Z} \times \mathbb{R}^{n-k}$ à support compact et $\psi_{k} \circ \tilde{D} \circ \psi_{k}^{-1}$ appartient au voisinage $\Omega$ de l'identité dans $Diff_{0}^{\infty}(V)$.

Lorsque $f$ est suffisamment proche de l'identité, $g(x,y)$ appartient, pour tous $x$ et $y$, au voisinage $W$ de l'identité donné par le lemme précédent.
On obtient alors que, pour $(x,y) \in \mathbb{R}^{k-1} \times \mathbb{R}^{n-k}$ :
$$g(x,y)=[a(g(x,y)),A][b(g(x,y)),B] \circ [c(g(x,y)),C]$$
puis que
$$g(x,y)=[a(g(x,y)),\tilde{A}(x,y)][b(g(x,y)),\tilde{B}(x,y)] \circ [c(g(x,y)),\tilde{C}(x,y)]$$
En effet, sur $[-1,1]^{n-1}$, on a $\tilde{A}(x,y)=A$, $\tilde{B}(x,y)=B$ et $\tilde{C}(x,y)=C$, et, hors de $[-1,1]^{n-1}$, les deux membres de l'égalité valent l'identité. Ainsi, le difféomorphisme $\psi_{k}^{-1} \circ \phi_{k}(f) \circ \psi_{k}$ s'écrit comme produit de 3 commutateurs d'éléments de $\Omega$, lorsque $f$ est suffisamment proche de l'identité. Un difféomorphisme $f$ proche de l'identité s'écrit alors comme produit de $3n$ commutateurs d'éléments de $\Omega$, ce qui démontre le théorème.

\end{document}